\documentclass{article}
\usepackage{amssymb}	% Allows use of AMS's mathematical symbols
\usepackage{amsmath}    % Allows use of AMS's mathematical commands
\usepackage{setspace} % use \doublespacing after \maketitle
\usepackage{latexsym}	% Allows use of old latex symbols
\usepackage{longtable}
\usepackage{multirow}
\usepackage{epsfig}
\newtheorem{definition}{Definition}[section]

\newtheorem{theorem}{Theorem}[section]

\newtheorem{remark}{Remark}[section]

\newcommand{\qed}{\nobreak \ifvmode \relax \else \ifdim\lastskip<1.5em \hskip-\lastskip \hskip1.5em plus0em minus0.5em \fi \nobreak \vrule height0.75em width0.5em depth0.25em\fi} 

\def\A{{\bf A}}
\def\B{{\bf B}}

\def\I{{\bf I}}
\def\J{{\bf J}}
\def\K{{\bf K}}

\def\0{{\bf 0}}

\def\b{{\bf b}}

\def\e{{\bf e}}

\def\g{{\bf g}}

\def\m{{\bf m}}

\def\q{{\bf q}}

\def\u{{\bf u}}

\def\x{{\bf x}}

\def\Tr{{\rm T}}

\title{Controllability of Spacecraft Using Only Magnetic Torques}
\author {Yaguang Yang\thanks{
Office of Research, NRC, 21 Church Street, Rockville, 20850. Email:
yaguang.yang@verizon.net} 
}
\date{\today}

\begin{document}

\maketitle  
%\doublespacing 

\begin{abstract}
%It is well-known that the magnetic torque coils interacting with 
%the Earth magnetic field 
%can provide torques in only two dimensional space at any moment. 
%But due to the change of the Earth magnetic field
%at the different time (or location) of the spacecraft, people 
%believe that the spacecraft attitude is still
%controllable by using only magnetic torque coils. Many designs 
%were proposed under the assumption
%of the controllability. Some serious efforts were made to 
%establish the conditions for the controllability to
%hold. But these conditions are either for a weak ``average'' 
%controllability or are not easy to be verified.
%In this paper, we will establish conditions under which the 
%conventional controllability holds for the 
%spacecraft using only magnetic torques, 
Spacecraft attitude control using only magnetic torques is a 
time-varying system. Many designs were proposed using LQR and 
$H_{\infty}$ formulations. The existence of the solutions depends 
on the controllability of the linear time-varying systems which 
has not been established. In this paper, we will derive the 
conditions of the controllability for this linear time-varying 
systems.
\end{abstract}


\newpage

\section{ Introduction}

Spacecraft attitude control using magnetic torque is a very attractive 
technique because the implementation is simple, the system is reliable (without moving 
mechanical parts), the torque coils are inexpensive, and their weights are light.
The main issue of using only magnetic torques to control the
attitude is that the magnetic torques generated by magnetic coils
are not available in all desired axes at any time \cite{sidi97}. However, because of the
constant change of the Earth's magnetic field as a spacecraft circles around
the earth, the controllable subspace changes all the time, many researchers 
believe that spacecraft's attitude is actually controllable by using only
magnetic torques. Numerous spacecraft attitude control designs were
proposed in the last twenty five years exploring the features of the time-varying
systems \cite{me89,pitte93,wisni97,psiaki01,la04,sl05,la06,yra07,plv10,cw10,rh11,zl11}.
Some of these papers tried Euler angle model and Linear Quadratic Regulator (LQR) formulations
\cite{me89,pitte93,wisni97,psiaki01,cw10} which are explicitly or implicitly assumed 
that the controllability for the linear time-varying system holds so that the optimal 
solutions exist \cite{kalman60}. But for the problem of spacecraft 
attitude control using only magnetic torque, no controllability condition has been established for this linear time-varying system
to the best of our knowledge. 

Other researchers \cite{la04,sl05,la06} proposed direct design 
methods using Lyapunov stabilization theory. The 
existence of the solutions for these methods implicitly 
depends on the controllability for the nonlinear time-varying 
system. Therefore, Bhat \cite{bhat05} investigated controllability 
of the nonlinear time-varying systems. However, the condition for
the controllability of the nonlinear time-varying systems
established in this paper  is hard to be verified and is a sufficient condition.

Recently, a reduced quaternion model was proposed in \cite{yang10} and its
 merits over Euler angle model were discussed in \cite{yang10,yang12,yang14}.
The reduced quaternion model was also used for the design of spacecraft attitude control 
system using magnetic torque \cite{plv10,rh11,zl11}. Because the
controllability of the linear time-varying systems was not established, 
the existence of the solutions was not guaranteed.

In this paper, we will consider the reduced linear quaternion model
proposed in \cite{yang10} and incorporate control model using only
magnetic torques. We will establish the conditions of the controllability for this 
very general linear time-varying system. The same strategy can easily
be used to prove the controllability of the Euler angle based linear time-varying
system considered in \cite{psiaki01}. However, we will not derive the similar result because of the merits of
the reduced quaternion model as discussed in \cite{yang10,yang12,yang14}.

The remainder of the paper is organized as follows. Section 2 
provides a description of the linear time-varying model of the
spacecraft attitude control system using only magnetic torque. 
Section 3 gives the proof of the controllability for this linear 
time-varying system. The conclusions are summarized in Section 4.

\section{The linear time-varying model} % for attitude control using only magnetic torque}

Let $\J$ be the inertia matrix of a spacecraft defined by
\begin{eqnarray}
\J =\left[   \begin{array}{ccc}
J_{11} & J_{12} & J_{13} \\
J_{21} & J_{22} & J_{23} \\
J_{31} & J_{32} & J_{33}   
\end{array} \right],
\label{inertia}
\end{eqnarray}
We will consider the nadir pointing spacecraft. Therefore, the attitude of the spacecraft is represented by the 
rotation of the spacecraft body frame relative to the local vertical and local horizontal 
(LVLH) frame. Therefore, we will represent the quaternion and spacecraft body rate in 
terms of the rotations of the spacecraft body frame relative to the LVLH frame. Let 
$\boldsymbol{\omega}=[\omega_1,\omega_2, \omega_3]^{\Tr}$ be the 
body rate with respect to the LVLH frame represented in the body frame, 
$\omega_0$ be the orbit (and LVLH frame) rate with respect to 
the inertial frame, represented in the LVLH frame.
Let $\bar{\q}=[q_0, q_1, q_2, q_3]^{\Tr}=[q_0, \q^{\Tr}]^{\Tr}=
[\cos(\frac{\alpha}{2}), \hat{\e}^{\Tr}\sin(\frac{\alpha}{2})]^{\Tr}$
be the quaternion representing the rotation of the body frame relative to the LVLH frame, 
where $\hat{\e}$ is the unit length rotational axis and $\alpha$ is the rotation angle about $\hat{\e}$. 
Therefore, the reduced kinematics equation becomes \cite{yang10}
\begin{eqnarray} \nonumber
\left[  \begin{array} {c} \dot{q}_1 \\ \dot{q}_2 \\ \dot{q}_3
\end{array} \right]
& = &\frac{1}{2}  \left[  \begin{array} {ccc} 
\sqrt{1-q_1^2-q_2^2-q_3^2} & -q_3  & q_2 \\
q_3 & \sqrt{1-q_1^2-q_2^2-q_3^2} & -q_1 \\
-q_2  &  q_1  & \sqrt{1-q_1^2-q_2^2-q_3^2} \\
\end{array} \right] 
\left[  \begin{array} {c}  \omega_{1} \\ \omega_{2} \\ \omega_{3}
\end{array} \right]  \\
& = & \g(q_1,q_2, q_3, \boldsymbol{\omega}).
\label{nadirModel2}
\end{eqnarray}
Assume that the inertia matrix of the spacecraft is diagonal which
is approximately correct for real systems, let the control 
torque vector be $\u=[u_x,u_y,u_z]^{\Tr}$, then the linearized 
nadir pointing spacecraft model with gravity gradient disturbance torque 
is given as follows \cite{yang10}:
\begin{eqnarray}
\left[  \begin{array}{c} 
\dot{q}_1 \\
\dot{q}_2 \\
\dot{q}_3 \\
\dot{\omega}_1 \\
\dot{\omega}_2 \\ 
\dot{\omega}_3 
\end{array}  \right] 
=\left[  \begin{array}{cccccc} 
0 & 0 & 0 & .5 & 0 & 0 \\
0 & 0 & 0 & 0 & .5 & 0 \\
0 & 0 & 0 & 0 & 0 & .5 \\
f_{41} & 0 & 0 & 0 & 0 & f_{46} \\
0 & f_{52} & 0 & 0 & 0 & 0 \\ 
0 & 0 & f_{63} & f_{64} & 0 & 0 
\end{array}  \right] 
\left[  \begin{array}{c} 
{q}_1 \\
{q}_2 \\
{q}_3 \\
{\omega}_1 \\
{\omega}_2 \\ 
{\omega}_3 
\end{array}  \right]
+\left[  \begin{array}{c} 
0 \\
0 \\
0 \\
u_x/J_{11} \\
u_y/J_{22} \\ 
u_z/J_{33} 
\end{array}  \right] \nonumber \\
\label{23}
\end{eqnarray}
where 
\begin{subequations}
\begin{gather}
f_{41}=[8(J_{33}-J_{22})\omega_0^2]/J_{11} \\
f_{46}=(-J_{11}+J_{22}-J_{33})\omega_0/J_{11} \\
f_{64}=(J_{11}-J_{22}+J_{33})\omega_0/J_{33} \\
f_{52}=[6(J_{33}-J_{11})\omega_0^2]/J_{22} \\
f_{63}=[2(J_{11}-J_{22})\omega_0^2]/J_{33}.
\end{gather}
\label{para}
\end{subequations}
The control torques generated by magnetic coils interacting with
the Earth's magnetic field is given by (see \cite{sidi97})
\[
\u=\m \times \b
\]
where the Earth's magnetic field in spacecraft coordinates, 
$\b(t)=[b_1(t),b_2(t),b_3(t)]^{\Tr}$, is computed using the 
spacecraft position, the spacecraft attitude, and a spherical harmonic model of the Earth's magnetic field \cite{wertz78}; and 
$\m=[m_1,m_2,m_3]^{\Tr}$ is the spacecraft magnetic coils' induced 
magnetic moment in the spacecraft coordinates.

The time-variation of the system is an approximate periodic function of $\b(t)=\b(t+T)$ where $T=\frac{2\pi}{\omega_0}$ is
the orbital period. This magnetic field $\b(t)$ can be 
approximately expressed as follows \cite{psiaki01}:
\begin{equation}
\left[  \begin{array}{c} 
b_1(t) \\
b_2(t) \\
b_3(t) 
\end{array}  \right] 
= \frac{\mu_f}{a^{3}}
\left[  \begin{array}{c} 
\cos(\omega_0t)\sin(i_m) \\
-\cos(i_m) \\
2\sin(\omega_0t)\sin(i_m) 
\end{array}  \right],
\label{field}
\end{equation}
where $i_m$ is the inclination of the spacecraft orbit with 
respect to the magnetic equator, $\mu_f=7.9 \times 10^{15}$
Wb-m is the field's dipole strength, and $a$ is the orbit's 
semi-major axis. The time $t=0$ is measured at the ascending-node
crossing of the magnetic equator. Therefore, the reduced
quaternion linear time-varying system is as follows:
\begin{eqnarray}
\left[  \begin{array}{c} 
\dot{q}_1 \\
\dot{q}_2 \\
\dot{q}_3 \\
\dot{\omega}_1 \\
\dot{\omega}_2 \\ 
\dot{\omega}_3 
\end{array}  \right] 
&=& \left[  \begin{array}{cccccc} 
0 & 0 & 0 & .5 & 0 & 0 \\
0 & 0 & 0 & 0 & .5 & 0 \\
0 & 0 & 0 & 0 & 0 & .5 \\
f_{41} & 0 & 0 & 0 & 0 & f_{46} \\
0 & f_{52} & 0 & 0 & 0 & 0 \\ 
0 & 0 & f_{63} & f_{64} & 0 & 0 
\end{array}  \right] 
\left[  \begin{array}{c} 
{q}_1 \\
{q}_2 \\
{q}_3 \\
{\omega}_1 \\
{\omega}_2 \\ 
{\omega}_3 
\end{array}  \right]
+\left[  \begin{array}{ccc} 
0 & 0 & 0 \\
0 & 0 & 0  \\
0 & 0 & 0  \\
0 & \frac{b_3(t)}{J_{11}} & -\frac{b_2(t)}{J_{11}} \\
-\frac{b_3(t)}{J_{22}} & 0 & \frac{b_1(t)}{J_{22}} \\ 
\frac{b_2(t)}{J_{33}} & -\frac{b_1(t)}{J_{33}} & 0
\end{array}  \right] 
\left[  \begin{array}{c} 
{m}_1 \\
{m}_2 \\
{m}_3 
\end{array}  \right] \nonumber \\
& := & 
\left[ \begin{array}{cc}
\0_3 & \frac{1}{2} \I_3  \\ \boldsymbol{\Lambda}_1 & \boldsymbol{\Sigma}_1
\end{array} \right]
\left[ \begin{array}{c}
\q \\ \boldsymbol{\omega} 
\end{array} \right]
+ \left[ \begin{array}{c}
\0_3 \\ \B_2(t)
\end{array} \right] \m
=\A\x+\B(t) \m.
\label{varying}
\end{eqnarray}
Substituting (\ref{field}) into (\ref{varying}) yields
\begin{equation}
\B_2(t)= \left[ \begin{array}{ccc}
0 & b_{42}(t) & b_{43}(t) \\
b_{51}(t) & 0 & b_{53}(t) \\ 
b_{61}(t) & b_{62}(t)  & 0
\end{array} \right]
\label{bt}
\end{equation}
where 
\begin{subequations}
\begin{gather}
b_{42} (t) = \frac{2\mu_f}{a^3 J_{11}} \sin(i_m) \sin(\omega_0t) 
\\
b_{43} (t)  = \frac{\mu_f}{a^3 J_{11}} \cos(i_m) 
\\
b_{53} (t)  = \frac{\mu_f}{a^3 J_{22}} \sin(i_m) \cos(\omega_0t) 
\\
b_{51} (t)  = -\frac{2\mu_f}{a^3 J_{22}} \sin(i_m) \sin(\omega_0t) 
= -b_{42}\frac{J_{11}}{J_{22}} 
\\
b_{61} (t) =-\frac{\mu_f}{a^3 J_{33}} \cos(i_m) = -b_{43} \frac{J_{11}}{J_{33}}
\\
b_{62} (t) =-\frac{\mu_f}{a^3 J_{33}} \sin(i_m) \cos(\omega_0t)
=-b_{53}\frac{J_{22}}{J_{33}}.
\end{gather}
\label{bij}
\end{subequations}
Therefore, we have
\begin{subequations}
\begin{gather}
b'_{42} (t)  = \frac{2\mu_f\omega_0}{a^3 J_{11}} \sin(i_m) \cos( \omega_0 t) 
\\
b'_{43} (t)  = 0
\\
b'_{53} (t)  = -\frac{\mu_f\omega_0}{a^3 J_{22}} \sin(i_m) \sin(\omega_0t) 
\\
b'_{51} (t)  = -\frac{2\mu_f\omega_0}{a^3 J_{22}} \sin(i_m) \cos(\omega_0t) 
= -b'_{42}\frac{J_{11}}{J_{22}} 
\\
b'_{61} (t) =0
\\
b'_{62} (t) =\frac{\mu_f\omega_0}{a^3 J_{33}} \sin(i_m) \sin(\omega_0t)
=-b_{53}\frac{J_{22}}{J_{33}}
\end{gather}
\label{bp}
\end{subequations}
and
\begin{subequations}
\begin{gather}
b''_{42} (t)  = -\frac{2\mu_f\omega_0^2}{a^3 J_{11}} \sin(i_m) \sin(\omega_0t) 
\\
b''_{43} (t)  = 0 
\\
b''_{53} (t)  = -\frac{\mu_f\omega_0^2}{a^3 J_{22}} \sin(i_m) \cos(\omega_0t) 
\\
b''_{51} (t)  =  \frac{2\mu_f\omega_0^2}{a^3 J_{22}} \sin(i_m) \sin(\omega_0t) 
= -b''_{42}\frac{J_{11}}{J_{22}} 
\\
b''_{61} (t) =0
\\
b''_{62} (t) = \frac{\mu_f\omega_0^2}{a^3 J_{33}} \sin(i_m) \cos(\omega_0t)
=-b''_{53}\frac{J_{22}}{J_{33}}.
\end{gather}
\label{bpp}
\end{subequations}
In matrix format, we have
\begin{equation}
\B'_2(t)= \left[ \begin{array}{ccc}
0 & b'_{42} & 0 \\
b'_{51} & 0 & b'_{53} \\ 
0 & b'_{62}  & 0
\end{array} \right],
\end{equation}
and
\begin{equation}
\B''_2(t)= \left[ \begin{array}{ccc}
0 & b''_{42} & 0 \\
b''_{51} & 0 & b''_{53} \\ 
0 & b''_{62}  & 0
\end{array} \right].
\end{equation}

A special case is when $i_m=0$, i.e., the spacecraft orbit is on
the equator plane of the Earth's magnetic field. In this case,
$\b(t)=[0,-\frac{\mu_f}{a^3},0]^{\Tr}$ is a constant vector. The 
linear time-varying system of this special case is reduced to a 
linear time-invariant system whose model is given by
\begin{eqnarray}
\left[  \begin{array}{c} 
\dot{q}_1 \\
\dot{q}_2 \\
\dot{q}_3 \\
\dot{\omega}_1 \\
\dot{\omega}_2 \\ 
\dot{\omega}_3 
\end{array}  \right] 
&=& \left[  \begin{array}{cccccc} 
0 & 0 & 0 & .5 & 0 & 0 \\
0 & 0 & 0 & 0 & .5 & 0 \\
0 & 0 & 0 & 0 & 0 & .5 \\
f_{41} & 0 & 0 & 0 & 0 & f_{46} \\
0 & f_{52} & 0 & 0 & 0 & 0 \\ 
0 & 0 & f_{63} & f_{64} & 0 & 0 
\end{array}  \right] 
\left[  \begin{array}{c} 
{q}_1 \\
{q}_2 \\
{q}_3 \\
{\omega}_1 \\
{\omega}_2 \\ 
{\omega}_3 
\end{array}  \right]
+\left[  \begin{array}{ccc} 
0 & 0 & 0 \\
0 & 0 & 0  \\
0 & 0 & 0  \\
0 & 0 & -b_2/J_{11} \\
0 & 0 & 0 \\ 
b_2/J_{33} & 0 & 0
\end{array}  \right] 
\left[  \begin{array}{c} 
{m}_1 \\
{m}_2 \\
{m}_3 
\end{array}  \right] \nonumber \\
& = & 
\A\x+\B \m.
\label{invariant}
\end{eqnarray}

\section{The proof of the controllability} % of the linear time-varying system}

The definition of controllability of linear
time-varying systems can be found in \cite[page 124]{rugh93}.
\begin{definition}
The linear state equation (\ref{varying}) is called controllable
on $[t_0,t_f]$ if given any $x_0$, there exists a continuous input
signal $\m(t)$ defined on $[t_0,t_f]$ such that the correspnding
solution of (\ref{varying}) satisfies $x(t_f)=0$.
\end{definition}
A main theorem used to prove the controllability of 
(\ref{varying}) is also given in \cite[page 127]{rugh93}.

\begin{theorem}
Let the state transition matrix $\boldsymbol{\Phi}(t,\tau)=e^{\A(t-\tau)}$.
Denote 
\begin{equation}
\K_j(t) = \frac{\partial^j}{\partial \tau^j} \left[ 
\boldsymbol{\Phi}(t, \tau)\B(\tau) \right] \Bigr|_{\tau = t}, 
\hspace{0.1in} j=1, 2, \ldots
\label{controllable}
\end{equation}
if $p$ is a positive integer such that, for $t \in [t_0,t_f]$, 
$\B(t)$ is $p$ time continuously differentiable. Then, the linear
time-varying equation (\ref{varying}) is controllable on 
$[t_0,t_f]$ if for some $t_c \in [t_0,t_f]$
\begin{equation}
\text{rank} \left[ \K_0(t_c), \K_1(t_c), \ldots, \K_p(t_c) \right]=n.
\label{rank}
\end{equation}
\label{theorem1}
\end{theorem}
\begin{remark}
If $\A$ and $\B$ are constant matrices, the rank condition of 
(\ref{rank}) for the linear time-varying system is reduced to the rank
condition for the linear time-invariant system \cite[page 128]{rugh93}, i.e.,
if
\begin{equation}
\text{rank} \left[ \B, \A\B, \ldots, \A^{n-1}\B\right]=n.
\label{rank1}
\end{equation}
then the linear time-invariant system $(\A, \B)$ is controllable.
\end{remark}
First, we consider the special case of (\ref{invariant}), 
the time-invariant system when the spacecraft
orbit is on the equator plane of the Earth's magnetic field ($i_m =0$). Let
$\boldsymbol{\Sigma}$ denote any $3\times 3$ diagonal or anti-diagonal matrix 
with the second row composed of zeros
\[
\boldsymbol{\Sigma} := \Bigg\{ \left[  \begin{array}{ccc} 
0 & 0 & \times  \\
0 & 0 & 0  \\
\times  & 0 & 0  
\end{array}  \right] \hspace{0.2in} 
\mbox{or} \hspace{0.2in} 
\left[  \begin{array}{ccc} 
\times  & 0 & 0 \\
0 & 0 & 0  \\
0 & 0 & \times   
\end{array}  \right] \Bigg\},
\]
and $\boldsymbol{\Lambda}$ denote any $3\times 3$ diagonal matrix with the form
\[
\boldsymbol{\Lambda} := \Bigg\{ \left[  \begin{array}{ccc} 
\times & 0 & 0 \\
0 & \times & 0  \\
0 & 0 & \times  
\end{array}  \right]  \Bigg\}.
\]
It is easy to verify that if 
$\boldsymbol{\Sigma}_i \in \boldsymbol{\Sigma}$, 
$\boldsymbol{\Sigma}_j \in \boldsymbol{\Sigma}$, 
and $\boldsymbol{\Lambda}_k \in \boldsymbol{\Lambda}$, 
then $\boldsymbol{\Sigma}_i \boldsymbol{\Sigma}_j \in \boldsymbol{\Sigma}$,
$\boldsymbol{\Sigma}_i+  \boldsymbol{\Sigma}_j \in \boldsymbol{\Sigma}$, and 
$\boldsymbol{\Lambda}_k \boldsymbol{\Sigma}_i \in 
\boldsymbol{\Sigma}$. Using this fact to expand
the matrix $[\B, \A^2\B, \A^3\B, \A^4\B, \A^5\B]$,
where $\A$ and $\B$ are defined in (\ref{invariant}), shows that
the second row of the controllability matrix in (\ref{rank1}) is composed of 
all zeros. 
This proves that if the spacecraft orbit is on the equator
plane of the Earth's magnetic field, the spacecraft attitude
cannot be stabilized by using only magnetic torques.

Now we show that under some simple conditions, the linear time-varying system 
(\ref{varying}) is controllable for any orbit which is not on the 
equator plane of the Earth's magnetic field, i.e., $i_m \neq 0$. From (\ref{controllable}), we have
\[
\K_0(t)=\boldsymbol{\Phi}(t,t)\B(t)=e^{\A(t-t)}\B(t)=\B(t),
\]
\begin{eqnarray}
\K_1(t)& = & \frac{\partial}{\partial \tau} \left[ 
\boldsymbol{\Phi}(t, \tau)\B(\tau) \right] \Bigr|_{\tau = t}
= \frac{\partial}{\partial \tau} \left[ 
e^{\A(t-\tau)}\B(\tau) \right] \Bigr|_{\tau = t} \nonumber \\
& = & \left[ -\A e^{\A(t-\tau )}\B(\tau ) + e^{\A(t-\tau )}\B'(\tau ) 
\right] \Bigr|_{\tau = t} = -\A\B(t)+\B'(t),  
\label{k1}
\end{eqnarray}
\begin{eqnarray}
\K_2(t)& = & \frac{\partial^2}{\partial \tau^2} \left[ 
\boldsymbol{\Phi}(t, \tau)\B(\tau) \right] \Bigr|_{\tau = t}
\nonumber \\
& = & \left[ 
\A^2e^{\A(t-\tau)}\B(\tau) -2\A e^{\A(t-\tau )}\B'(\tau ) + e^{\A(t-\tau )}\B''(\tau ) 
\right] \Bigr|_{\tau = t} \nonumber \\
& = & \A^2\B(t) -2\A\B'(t)+\B''(t). 
\label{k2} 
\end{eqnarray}

Using the notation of (\ref{varying}), we can rewrite equation (\ref{k1}) as 
\[
\K_1(t) = -\left[ \begin{array}{cc}
\0_3 & \frac{1}{2} \I_3  \\ \boldsymbol{\Lambda}_1 & \boldsymbol{\Sigma}_1
\end{array} \right]
\left[ \begin{array}{c}
\0_3  \\ \B_2 
\end{array} \right]+
\left[ \begin{array}{c}
\0_3  \\ \B'_2 
\end{array} \right]
=\left[ \begin{array}{c}
-\frac{1}{2}\B_2 \\ -\boldsymbol{\Sigma}_1 \B_2 +\B'_2
\end{array} \right].
\]
Since
\[
\A^2\B= \A \left[ \begin{array}{cc}
\0_3 & \frac{1}{2} \I_3  \\ \boldsymbol{\Lambda}_1 & \boldsymbol{\Sigma}_1
\end{array} \right]
\left[ \begin{array}{c}
\0_3  \\ \B_2 
\end{array} \right]
=\left[ \begin{array}{cc}
\0_3 & \frac{1}{2} \I_3  \\ \boldsymbol{\Lambda}_1 & \boldsymbol{\Sigma}_1
\end{array} \right]
\left[ \begin{array}{c}
\frac{1}{2}\B_2 \\ \boldsymbol{\Sigma}_1 \B_2 
\end{array} \right]
=\left[ \begin{array}{c}
\frac{1}{2}\boldsymbol{\Sigma}_1 \B_2 \\ \frac{1}{2} \boldsymbol{\Lambda}_1 \B_2+\boldsymbol{\Sigma}_1^2 \B_2 
\end{array} \right]
\]
and
\[
-2\A \B'=-2 \left[ \begin{array}{cc}
\0_3 & \frac{1}{2} \I_3  \\ \boldsymbol{\Lambda}_1 & \boldsymbol{\Sigma}_1
\end{array} \right] \left[ \begin{array}{c}
\0_3  \\ \B'_2 
\end{array} \right]
= \left[ \begin{array}{c}
-\B'_2 \\ -2 \boldsymbol{\Sigma}_1 \B'_2
\end{array} \right],
\]
equation (\ref{k2}) is reduced to 
\[
\K_2(t) = \A^2\B-2\A\B'+\B''= \left[ \begin{array}{c}
\frac{1}{2}\boldsymbol{\Sigma}_1 \B_2 -\B'_2  \\ \frac{1}{2} \boldsymbol{\Lambda}_1 \B_2+\boldsymbol{\Sigma}_1^2 \B_2-2 \boldsymbol{\Sigma}_1 \B'_2+\B''_2
\end{array} \right].
\]
Hence,
\begin{eqnarray}
[\K_0(t), \K_1(t), \K_2(t)] & = &  [ \B(t) \,\,\, |  \,\,\, -\A\B(t)+\B'(t) \,\,\, |  \,\,\, \A^2\B(t)-2\A\B'(t)+\B''(t)] 
\nonumber \\
& = & \left[ 
\begin{smallmatrix}
\0_3  \\  \B_2  
\end{smallmatrix} \middle|  
\begin{smallmatrix}
 -\frac{1}{2} \B_2  \\ -\boldsymbol{\Sigma}_1 \B_2+\B'_2  
 \end{smallmatrix}  \middle|  
\begin{smallmatrix} 
\frac{1}{2}\boldsymbol{\Sigma}_1 \B_2 -\B'_2(t)  \\  \frac{1}{2} \boldsymbol{\Lambda}_1 \B_2+\boldsymbol{\Sigma}_1^2 \B_2 -2\boldsymbol{\Sigma}_1 \B'_2 +\B''_2 
\end{smallmatrix}
\right].
\end{eqnarray}
Notice that
\begin{eqnarray}
& &  \text{rank} [\K_0(t), \K_1(t), \K_2(t)] \nonumber \\
& = & \text{rank}  \left( 
\left[ \begin{array}{cc}
\I_3 & 0_3    \\ -2 \boldsymbol{\Sigma}_1 & \I_3
\end{array} \right]
\left[ 
\begin{smallmatrix}
\0_3  \\  \B_2  
\end{smallmatrix} \middle|  
\begin{smallmatrix}
 -\frac{1}{2} \B_2  \\ -\boldsymbol{\Sigma}_1 \B_2+\B'_2  
 \end{smallmatrix}  \middle|  
\begin{smallmatrix} 
\frac{1}{2}\boldsymbol{\Sigma}_1 \B_2 -\B'_2(t)  \\  \frac{1}{2} \boldsymbol{\Lambda}_1 \B_2+\boldsymbol{\Sigma}_1^2 \B_2 -2\boldsymbol{\Sigma}_1 \B'_2 +\B''_2 
\end{smallmatrix}
\right] \right)  \nonumber \\
& = &  \text{rank} \left[ 
\begin{smallmatrix}
\0_3  \\  \B_2  
\end{smallmatrix} \middle|  
\begin{smallmatrix}
 -\frac{1}{2} \B_2  \\ \B'_2  
 \end{smallmatrix}  \middle|  
\begin{smallmatrix} 
\frac{1}{2}\boldsymbol{\Sigma}_1 \B_2 -\B'_2(t)  \\  \frac{1}{2} \boldsymbol{\Lambda}_1 \B_2 +\B''_2 
\end{smallmatrix}
\right] \nonumber \\
& = &  \text{rank} \left[ 
\begin{smallmatrix}
\0_3  \\  \B_2  
\end{smallmatrix} \middle|  
\begin{smallmatrix}
 - \B_2  \\ \B'_2  
 \end{smallmatrix}  \middle|  
\begin{smallmatrix} 
\boldsymbol{\Sigma}_1 \B_2 -2\B'_2(t)  \\  \frac{1}{2} \boldsymbol{\Lambda}_1 \B_2 +\B''_2 
\end{smallmatrix}
\right],
\end{eqnarray}

\begin{eqnarray}
& & \boldsymbol{\Sigma}_1 \B_2-2 \B'_2(t)
\nonumber \\
& = & \left[ \begin{array}{ccc}
0 & 0 & f_{46} \\
0 & 0 & 0 \\
f_{64} & 0 & 0
\end{array} \right]
\left[ \begin{array}{ccc}
0 & b_{42}(t) & b_{43}(t) \\
b_{51}(t) & 0 & b_{53}(t) \\ 
b_{61}(t) & b_{62}(t)  & 0
\end{array} \right]
-2\left[ \begin{array}{ccc}
0 & b'_{42} & 0 \\
b'_{51} & 0 & b'_{53} \\ 
0 & b'_{62}  & 0
\end{array} \right] \nonumber \\
& = & \left[ \begin{array}{ccc}
f_{46} b_{61}(t) & f_{46} b_{62}(t)-2b'_{42} & 0  \\
-2 b'_{51}  & 0 & 2b'_{53} \\
0 &  f_{64} b_{42}(t)-2b'_{62} & f_{64} b_{43}(t) 
\end{array} \right],
\end{eqnarray}
and
\begin{eqnarray}
& & \frac{1}{2} \boldsymbol{\Lambda}_1 \B_2 + \B''_2(t)
\nonumber \\
& = & \frac{1}{2} \left[ \begin{array}{ccc}
f_{41} & 0 & 0 \\
0 & f_{52}  & 0 \\
0 & 0 & f_{63} 
\end{array} \right]
\left[ \begin{array}{ccc}
0 & b_{42}(t) & b_{43}(t) \\
b_{51}(t) & 0 & b_{53}(t) \\ 
b_{61}(t) & b_{62}(t)  & 0
\end{array} \right]
+\left[ \begin{array}{ccc}
0 & b''_{42} & 0 \\
b''_{51} & 0 & b''_{53} \\ 
0 & b''_{62}  & 0
\end{array} \right] \nonumber \\
& = & \left[ \begin{array}{ccc}
0 & \frac{1}{2} f_{41} b_{42}(t) + b''_{42} & \frac{1}{2} f_{41} b_{43}(t)   \\
\frac{1}{2} f_{52} b_{51}(t) + b''_{51}  & 0 & \frac{1}{2} f_{52} b_{53}(t) + b''_{53} \\
\frac{1}{2} f_{63} b_{61}(t) &  \frac{1}{2} f_{63} b_{62}(t) + b''_{62} & 0
\end{array} \right],
\end{eqnarray}
we have
\begin{eqnarray}
& &  \left[ \begin{smallmatrix}
\0_3  \\  \B_2  
\end{smallmatrix} \middle|  
\begin{smallmatrix}
 - \B_2  \\ \B'_2  
 \end{smallmatrix}  \middle|  
\begin{smallmatrix} 
\boldsymbol{\Sigma}_1 \B_2 -2\B'_2(t)  \\  \frac{1}{2} \boldsymbol{\Lambda}_1 \B_2 +\B''_2 
\end{smallmatrix}
\right] \nonumber \\
& = &   
\footnotesize
\left[ \begin{array}{ccccccccc}
0 & 0 & 0 & 0 & -b_{42}(t) & -b_{43}(t) & f_{46} b_{61}(t) & f_{46} b_{62}(t)-2b'_{42} & 0  \\
0 & 0 & 0 & -b_{51}(t) & 0 & -b_{53}(t) & -2 b'_{51}  & 0 & 2b'_{53} \\
0 & 0 & 0 & -b_{61}(t) & -b_{62}(t)  & 0 & 0 &  f_{64} b_{42}(t)-2b'_{62} & f_{64} b_{43}(t) \\
0 & b_{42}(t) & b_{43}(t) & 0 & b'_{42} & 0  & 0  & \frac{1}{2} f_{41} b_{42}(t) + b''_{42} & \frac{1}{2} f_{41} b_{43}(t) \\
b_{51}(t) & 0 & b_{53}(t) & b'_{51} & 0 & -b_{53}(t) & \frac{1}{2} f_{52} b_{51}(t) + b''_{51}  & 0 & \frac{1}{2} f_{52} b_{53}(t) + b''_{53} \\
b_{61}(t) & b_{62}(t) & 0 & 0 & b'_{62}  & 0 & \frac{1}{2} f_{63} b_{61}(t) &  \frac{1}{2} f_{63} b_{62}(t) + b''_{62} & 0
\end{array} \right]. \nonumber
\normalsize
\end{eqnarray}
To show that this matrix is full rank for some $t_c$, 
we show that there is a $6 \times 6$ submatrix whose determinant 
is not zero for $\omega_0t_c = \frac{\pi}{2}$. In view of 
(\ref{bij}), (\ref{bp}), and (\ref{bpp}), for this $t_c$, 
we have
\begin{equation}
b_{53}(t_c)=b_{62}(t_c)=b'_{51}(t_c)=b'_{42}(t_c)=b''_{53}(t_c)=b''_{62}(t_c)=0.
\label{reduce}
\end{equation}
Consider the submatrix composed of the $1$st, $2$nd, $4$th, 
$5$th, $7$th, $8$th columns, and using (\ref{reduce}), we have 
\begin{eqnarray}
& & \det \left[ \begin{array}{cccccc}
0 & 0  & 0 & -b_{42}(t_c) & f_{46} b_{61}(t_c) & f_{46} b_{62}(t_c)-2b'_{42} \\
0 &  0 & -b_{51}(t_c) & 0 & -2 b'_{51}  & 0 \\
0 & 0 & -b_{61}(t_c)  & -b_{62}(t_c) & 0 &  f_{64} b_{42}(t_c)-2b'_{62} \\
0 & b_{42}(t_c)  & 0 & b'_{42} & 0  & \frac{1}{2} f_{41} b_{42}(t_c) + b''_{42} \\
b_{51}(t_c) & 0  & b'_{51} & 0 & \frac{1}{2} f_{52} b_{51}(t_c) + b''_{51}  & 0 \\
b_{61}(t_c) & b_{62}(t_c)  & 0 & b'_{62}  & \frac{1}{2} f_{63} b_{61}(t_c) &  \frac{1}{2} f_{63} b_{62}(t_c) + b''_{62} 
\end{array} \right] \nonumber \\
& = & \det \left[ \begin{array}{cccccc}
0 & 0  & 0 & -b_{42}(t_c) & f_{46} b_{61}(t_c) & 0 \\
0 &  0 & -b_{51}(t_c) & 0 & 0 & 0 \\
0 & 0 & -b_{61}(t_c)  & 0 & 0 &  f_{64} b_{42}(t_c)-2b'_{62} \\
0 & b_{42}(t_c)  & 0 & 0 & 0  & \frac{1}{2} f_{41} b_{42}(t_c) + b''_{42} \\
b_{51}(t_c) & 0  & 0& 0 & \frac{1}{2} f_{52} b_{51}(t_c) + b''_{51}  & 0 \\
b_{61}(t_c) & 0 & 0 & b'_{62}  & \frac{1}{2} f_{63} b_{61}(t_c) &  0
\end{array} \right] \nonumber \\
& = &  b_{51}(t_c) \det\left[ \begin{array}{ccccc}
0 & 0  & -b_{42}(t_c) & f_{46} b_{61}(t_c) & 0 \\
0 & 0 & 0 & 0 &  f_{64} b_{42}(t_c)-2b'_{62} \\
0 & b_{42}(t_c)  & 0 & 0  & \frac{1}{2} f_{41} b_{42}(t_c) + b''_{42} \\
b_{51}(t_c) & 0  & 0 & \frac{1}{2} f_{52} b_{51}(t_c) + b''_{51}  & 0 \\
b_{61}(t_c) & 0 & b'_{62}  & \frac{1}{2} f_{63} b_{61}(t_c) &  0
\end{array} \right] \nonumber \\
& = &  -\left( f_{64} b_{42}(t_c)-2b'_{62} \right) b_{51}(t_c) \det\left[ \begin{array}{cccc}
0 & 0  & -b_{42}(t_c) & f_{46} b_{61}(t_c) \\
0 & b_{42}(t_c)  & 0 & 0  \\
b_{51}(t_c) & 0  & 0 & \frac{1}{2} f_{52} b_{51}(t_c) + b''_{51}  \\
b_{61}(t_c) & 0 & b'_{62}  & \frac{1}{2} f_{63} b_{61}(t_c) 
\end{array} \right] \nonumber \\
& = &  - b_{42}(t_c)  \left( f_{64} b_{42}(t_c)-2b'_{62} \right) b_{51}(t_c) \det\left[ \begin{array}{ccc}
0 &  -b_{42}(t_c) & f_{46} b_{61}(t_c) \\
b_{51}(t_c) & 0 & \frac{1}{2} f_{52} b_{51}(t_c) + b''_{51}  \\
b_{61}(t_c) & b'_{62}  & \frac{1}{2} f_{63} b_{61}(t_c) 
\end{array} \right] \nonumber \\
& = & - b_{42}(t_c)  \left( f_{64} b_{42}(t_c)-2b'_{62} \right) 
b_{51}(t_c) \times \nonumber \\
& & \left[
b_{51} b'_{62} f_{46} b_{61}
-b_{42}\left( \frac{1}{2} f_{52}b_{51}+b''_{51} \right)b_{61}
+\frac{1}{2} f_{63}b_{61}b_{42}b_{51}
\right].
\end{eqnarray}
Therefore, in view of Theorem \ref{theorem1}, the time-varying system is controllable if 
\begin{equation}
f_{64} b_{42}(t_c)-2b'_{62} \neq 0,
\label{neq1}
\end{equation}
and
\begin{equation}
b_{51} b'_{62} f_{46} b_{61}
-b_{42}\left( \frac{1}{2} f_{52}b_{51}+b''_{51} \right)b_{61}
+\frac{1}{2} f_{63}b_{61}b_{42}b_{51} \neq 0.
\label{neq2}
\end{equation}
Using (\ref{para}), (\ref{bij}), (\ref{bp}), (\ref{bpp}), and 
noticing that $\sin(\omega_0 t_c) = \sin(\frac{\pi}{2})=1$,
we have
\begin{eqnarray}
& & f_{64} b_{42}(t_c)-2b'_{62} 
\nonumber \\
& = & \frac{(J_{11}-J_{22}+J_{33}) \omega_0}{J_{33}}
\frac{2\mu_f }{a^3 J_{11}}\sin(i_m) 
-2 \frac{\mu_f \omega_0}{a^3 J_{33}}\sin(i_m) \nonumber \\
& = & \frac{2\mu_f \omega_0 \sin(i_m)}{a^3 (J_{11}J_{33})} 
(J_{33}-J_{22}),  \nonumber
\end{eqnarray}
the first condition (\ref{neq1}) is reduced to
\begin{equation}
J_{33} \neq J_{22}.
\label{cond1}
\end{equation}
Repeatedly using the same relations, we have
\begin{eqnarray}
b_{51} b'_{62} f_{46} b_{61} & = & \left( -\frac{2\mu_f}{a^3 J_{22}} \sin(i_m) \right)
\left( \frac{\mu_f\omega_0}{a^3 J_{33}} \sin(i_m) \right) 
\left( \frac{(-J_{11}+J_{22}-J_{33})\omega_0}{J_{11}} \right)
\left(-\frac{\mu_f}{a^3 J_{33}} \cos(i_m) \right)
\nonumber \\
& = &  \frac{2 \mu_f^3 \omega_0^2 (-J_{11}+J_{22}-J_{33})}{a^9 J_{11}J_{22}J_{33}^2}
 \sin^2(i_m) \cos(i_m) ,
\label{1st}
\end{eqnarray}
\begin{eqnarray}
-b_{42}\left( \frac{1}{2} f_{52}b_{51}+b''_{51} \right)b_{61}
& = & - \left(  \frac{2\mu_f}{a^3 J_{11}} \sin(i_m) \right) 
\left( \frac{3(J_{33}-J_{11})\omega_0^2}{J_{22}}
\left( -\frac{2\mu_f}{a^3 J_{22}} \sin(i_m) \right) 
+\frac{2\mu_f\omega_0^2}{a^3 J_{22}} \sin(i_m)  \right) 
\nonumber \\
& & \left( -\frac{\mu_f}{a^3 J_{33}} \cos(i_m)  \right)
\nonumber \\
& = &   - \left(  \frac{2\mu_f}{a^3 J_{11}} \sin(i_m) \right) 
\left( \frac{2\mu_f\omega_0^2 }{a^3 J_{22}^2}\sin(i_m) ( -3J_{33} + 3J_{11} +J_{22}) \right)  
\nonumber \\
& & \left( -\frac{\mu_f}{a^3 J_{33}} \cos(i_m)  \right) \nonumber \\
& = &   \frac{4\mu_f^3 \omega_0^2 ( -3J_{33} + 3J_{11} +J_{22}) }{a^9  J_{11} J_{22}^2 J_{33}} \sin^2(i_m) \cos(i_m) ,
\label{2nd}
\end{eqnarray}
and
\begin{eqnarray}
\frac{1}{2} f_{63}b_{61}b_{42}b_{51}
& = &  \frac{ (J_{11}-J_{22})\omega_0^2}{J_{33}}
\left( -\frac{\mu_f}{a^3 J_{33}} \cos(i_m)  \right)
\left(  \frac{2\mu_f}{a^3 J_{11}} \sin(i_m) \right)
 \left(  -\frac{2\mu_f}{a^3 J_{22}} \sin(i_m) \right)
 \nonumber \\
& = &  \frac{4\mu_f^3 \omega_0^2 (J_{11}-J_{22}) }{ a^9  J_{11} J_{22} J_{33}^2} \sin^2(i_m) \cos(i_m).
\label{3rd}
\end{eqnarray}
Combining (\ref{1st}), (\ref{2nd}), and (\ref{3rd}), we can 
rewrite (\ref{neq2}) as
\begin{eqnarray}
& & b_{51} b'_{62} f_{46} b_{61}
-b_{42}\left( \frac{1}{2} f_{52}b_{51}+b''_{51} \right)b_{61}
+\frac{1}{2} f_{63}b_{61}b_{42}b_{51}  \nonumber \\
& = &   \frac{\mu_f^3 \omega_0^2}{ a^9  J_{11} J_{22}^2 J_{33}^2} \sin^2(i_m) \cos(i_m)
\nonumber \\
& & \left( 2J_{22}(-J_{11}+J_{22}-J_{33})+4J_{33} ( -3J_{33} + 3J_{11} +J_{22})+4 J_{22} (J_{11}-J_{22})  \right)
\nonumber \\
& = &   \frac{2\mu_f^3 \omega_0^2}{ a^9  J_{11} J_{22}^2 J_{33}^2} \sin^2(i_m) \cos(i_m)
[J_{22} (J_{11}-J_{22}+J_{33})-6J_{33} (J_{33}-J_{11})].
\label{all}
\end{eqnarray}
Therefore, the second condition of (\ref{neq2}) is reduced to
\begin{equation}
J_{22} (J_{11}-J_{22}+J_{33}) \neq 6J_{33} (J_{33}-J_{11}).
\end{equation}

We summarize our main result of this paper as the main 
theorem of this paper. 
\begin{theorem}
For the linear time-varying spacecraft attitude control system (\ref{varying}) using only magnetic torques,
if the orbit is on the equator plane of the Earth's magnetic field, then the spacecraft attitude is not fully 
controllable. If the orbit is not on the equator plane of the Earth's magnetic field, and the following two
conditions hold:
\begin{subequations}
\begin{gather}
J_{33} \neq J_{22}, \\
J_{22} (J_{11}-J_{22}+J_{33}) \neq 6J_{33} (J_{33}-J_{11}),
\end{gather}
\end{subequations}
then  the spacecraft attitude is fullly controllable by magnetic
coils. 
\end{theorem}

\begin{remark}
The controllability conditions include only the spacecraft orbit plane and the spacecraft 
inertia matrix which are very easy to be verified.
\end{remark}

\section{Conclusions} 
 
In this paper, the controllability of spacecraft attitude control 
system using only magnetic torques is considered. The conditions 
of the controllability are derived. These conditions 
are easier to be verifed than previously established conditions.

%\section{Acknowledgments} 

\end{document}